\documentclass[12pt,twoside]{amsart}

\usepackage{amsfonts}
\usepackage{amsmath}
\usepackage{epsfig}
\usepackage{latexsym}

\newtheorem{thm}{{\sc Theorem}}
\newtheorem{prop}{{\sc Proposition}}
\newtheorem{lem}{{\sc Lemma}}
\newtheorem{cor}{{\sc Corollary}}
\newtheorem{defn}{{\sc Definition}}

\newcommand{\qqed}{\hspace*{\fill} $\Box$}

\title{On singular fibres of complex Lagrangian fibrations}
\author{Daisuke Matsushita}
\address{Research Institue for Mathematical Sciences \\
         Kyoto University, Oiwake-Cyo Kitashirakawa \\
         Sakyo-Ku Kyoto 606-8052 Japan}
\thanks{*Research Fellow of the Japan Society for the Promotion of Science} 
\email{tyler@kurims.kyoto-u.ac.jp}

\begin{document}
\maketitle
\begin{abstract}
 We classify singular fibres over general points
 of the discriminant locus
 of projective 
 complex Lagrangian fibrations
 on 4-dimensional holomorphic symplectic manifolds. 
 The singular fibre $F$ is  the
 following either one: 
 $F$ is isomorphic to
 the product of an elliptic curve and a Kodaira singular fibre 
 up to finite unramified covering or 
 $F$ is a normal crossing variety consisting of
 several copies of a minimal elliptic ruled surface of which 
 the dual graph 
 is Dynkin diagram of type $A_n$, $\tilde{A_n}$  
 $\tilde{D_n}$ or $D_n$.
\end{abstract}
\section{Introduction} 
 First we define {\it complex Lagrangian fibration\/}.
\begin{defn}
 Let $(X, \omega )$ be a K\"{a}hler manifold with a
 holomorhpic symplectic two form $\omega$
 and $S$ a smooth manifold. A proper flat surjective morphism $f : X \to S$
 is said to be a complex Lagrangian fibration if 
 a general fibre $F$ of $f$ is a Lagrangian submanifold with
 respect to $\omega$,
 that is,
 the restriction of 2-form $\omega |_{F}$ is 
 identically zero and $\dim F = (1/2)\dim X$.
\end{defn}
\noindent
{\sc Remark. \quad} A general fibre $F$ of a complex
 Lagrangian fibration is a complex torus by Leauville's theorem. 

\vspace{5mm}

\begin{figure}[h]
    \epsffile{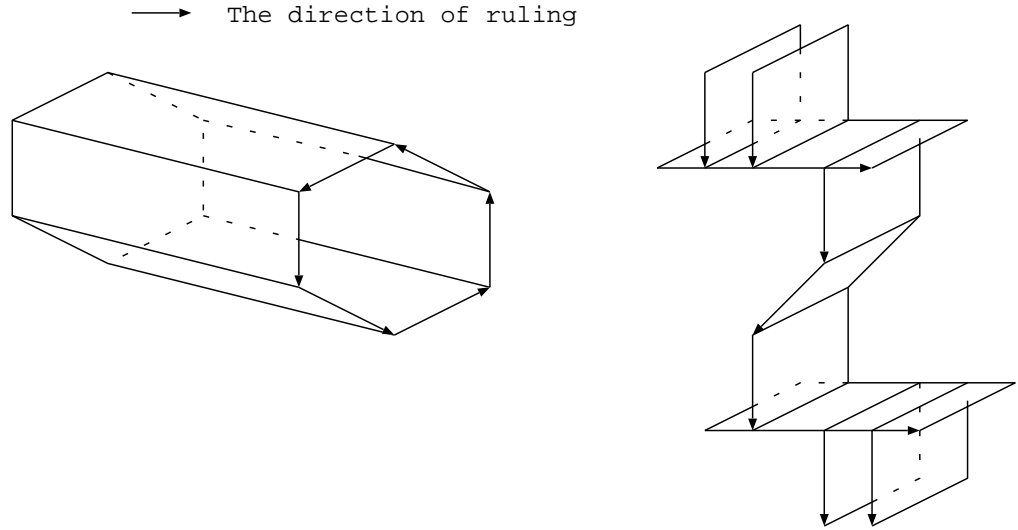}
    \caption{Figures of $\tilde{A_n}$ 
    and $\tilde{D_n}$ case.}
    \label{figure_tilde_A_n}
\end{figure}

\begin{figure}[h]
    \epsffile{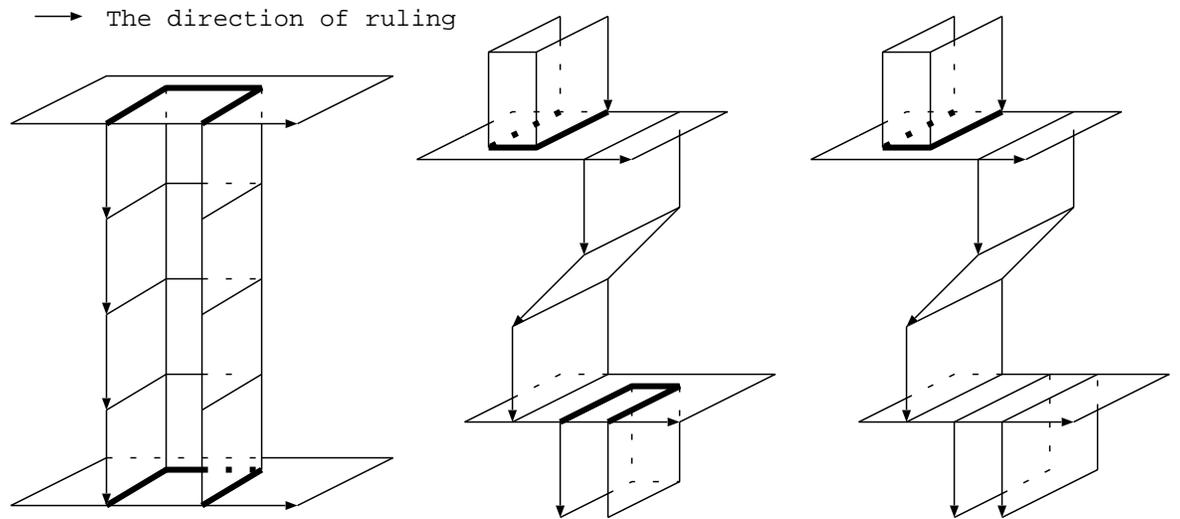}
    \caption{Figures of $A_n$, $A_n$ and $D_n$ case. Bold line
              represents bisection and line represents section.}
    \label{figure_A_n}
\end{figure}

\noindent
 The plainest example of a complex Lagrangian fibration
 is an elliptic fibration of  $K3$ surface over ${\mathbb P}^{1}$.
 In higher dimension,
 such a fibre space naturally occurs on a fiber space of
 an irreducible symplectic manifold ( \cite[Theorem 2]{matsu} 
 and \cite[Theorem 1]{matsu2} ). When the dimension of fibre
 is one, a complex Lagrangian fibration is 
 a minimal elliptic fibration and 
 whose  singular fibre 
 is completely 
 classified by Kodaira \cite[Theorem 6.2]{kodaira}.  
 In this note, we investigate
 singular fibres of a projective
 complex Lagrangian fibration whose fibre is 2-dimensional.
\begin{thm}\label{main}
 Let $f : X \to S$ be a complex Lagrangian fibration 
 on 4 dimensional symplectic manifold and
 $D$ the discriminant locus of $f$. Assume that $f$ is 
 projective.
 For a general point
 $x$ of $D$, 
 $f^{-1}(x)$ is the
 following one:
\begin{enumerate}
 \item There exists an \'{e}tale finite covering 
       $f^{-1}(x)^{\sim} \to f^{-1}(x)$ and 
       $f^{-1}(x)^{\sim}$ is
       isomorphic to the product of an elliptic curve and
       a Kodaira singular fibre of type 
       $I_0$, $I^{*}_{0}$, $II$, $II^{*}$, $III$, $III^{*}$,
       $IV$ or $IV^{*}$.
 \item $f^{-1}(x)$ is isomorphic to a
       normal crossing variety consisting of 
       a  minimal elliptic ruled surface. The
       dual graph of $f^{-1}(x)$
       is the Dynkin diagram of type $A_n$, $\tilde{A_n}$
       $D_n$ or $\tilde{D_n}$. 
       If the dual graph is of type $\tilde{A_n}$ or $\tilde{D_n}$,
       each double curve is 
       a section of the ruling.
       In the other cases, the double curve  
       on each edge components
       is a bisection or a section and  other double curve is 
       a section. 
       ( See figures \ref{figure_tilde_A_n} 
         and \ref{figure_A_n} (pp \pageref{figure_tilde_A_n}). )
\end{enumerate}
\end{thm}

\noindent Combining Theorem \ref{main} with
 \cite[Theorem 2]{matsu} and \cite[Theorem 1]{matsu2},
 we obtain the following corollary.
\begin{cor}
 Let $f :X \to B$ be a fibre space of a projective irreducible
 symplectic manifold. Assume that $\dim X = 4$. Then,
 for a general point $x$ of the discriminat locus of $f$,
 $f^{-1}(x)$ satisfies the properties of Theorem \ref{main} (1) or (2).
\end{cor}
\noindent
{\sc Remark. \quad} 
 Let $S$ be a $K3$ surface and $\pi : S \to {\mathbb P}^{1}$ 
 an elliptic fibration. The induced morhpism 
 $f : {\rm Hilb}^{2}S \to {\mathbb P}^{2}$ gives
 examples of singular fibres above except whose 
 dual graphs are $A_n$ or $D_n$.
 The author does not know whether a normal crossing variety
 whose dual graph is $A_n$ or $D_n$ occur as a 
 singular fibre
 of a fibre space of an irreducible symplectic manifold.

\vspace{5mm}

\noindent
 This paper is organized as follows. In section 2, we set up the 
 proof of Theorem \ref{main}. The key proposition is stated and
 proved in section 3. Section 4 and 5 are devoted to the proof
 of Theorem \ref{main}.

\vspace{5mm}
\noindent
{\sc Acknowledgment. \quad} 
 The author express his thanks to Professors A.~Fujiki, Y.~Miyaoka,
 S.~Mori 
 and  N.~Nakayama for their advice and encouragement.

\section{Preliminary}
\noindent
(2.1) In this section, we collect definitions and some fundamental
 materials which are necessary for the proof of Theorem \ref{main}. 
\begin{defn}
 Let $f : X \to \Delta^{1}$ be a proper surjective morphism from
 algebraic variety to an unit disk. 
 $f$ is said to be {\it semistable degeneration\/}
 if $f$ satisfies the following two properties:
\begin{enumerate}
 \item $f$ is smooth over $\Delta^{1}\setminus 0$.
 \item $f^{*}\{0 \}$ is a reduced normal crossing divisor.
\end{enumerate}
\end{defn}
\begin{defn}
 Let $f :X \to \Delta^{1}$ and 
     $f' : X' \to \Delta^{1}$ be proper surjective 
     morphisms from algebraic varieties to unit disks.
     We call $f$ is isomorphic to $f'$ (resp. birational)
     if there exists an isomorphism (resp. a birational map)
     $g : X \to X'$ such that $f' \circ g = f$.
\end{defn}
\noindent
(2.2) We refer the fundamental properties of an abelian fibration.
\begin{lem}\label{basic_properties}
 Let $f : (X,\omega ) \to S$ be a projective complex Lagrangian
 fibration.
\begin{enumerate}
 \item The discriminant locus of $f$ is pure codimension one.
 \item Let $F$ be an irreducible component of a fibre of $f$ and
       $j : \tilde{F} \to F$ a resolution of $F$.
       Then $j^{*}(\omega |_{F}) = 0$.
\end{enumerate}
\end{lem}
\noindent
{\sc Proof. \quad}
\begin{enumerate}
 \item Let $h : {\mathcal A} \to \Delta^{2}$ be a
 projective  flat morphism over a 
 two dimensional polydisk.
 Assume that a general fibre of $h$ is an abelian variety and
 $h$ is smooth over $\Delta^{2}\setminus 0$.
 For the proof of Lemma \ref{basic_properties} (1), it is enough to
 prove that $h$ is smooth morphism. 
 Let ${\mathcal A}^{\circ} := {\mathcal A}\setminus h^{-1}(0)$.
 Since $h$ is smooth 
 over $\Delta^{2}\setminus 0$ and 
 projective, there exists
 an \'{e}tale finite cover 
 $\pi : {\mathcal A}^{\circ} \to \bar{{\mathcal A}}^{\circ} $ and
 a smooth abelian fibration
 $\bar{h}^{\circ} : \bar{{\mathcal A}}^{\circ} \to \Delta^{2}\setminus 0$
 with a section.
 Since $\Delta^{2}\setminus 0$ is simply connected,
 $\bar{{\mathcal A}}^{\circ}$ has
 a level $n$ structure.
 We put ${\mathcal M}_{g}[n]$ the moduli space of
 $g$-dimensional abelian varieties with level $n$-structure. 
 There exists a morphism from
 $t : \Delta^{2}\setminus 0 \to {\mathcal M}_{g}[n]$
 and we extend $t$ on $\Delta^{2}$ by Hartogs theorem.
 The moduli ${\mathcal M}_{g}[n]$ has the universal family
 ${\mathcal A}_{g}[n] \to {\mathcal M}_{g}[n]$.
 Considering the pull back of ${\mathcal A}_{g}[n]$ by $t$,
 we obtain a smooth abelian morphism 
 $\bar{h} : \bar{{\mathcal A}} \to \Delta^{2}$ which is the extension of
 $\bar{h}^{\circ} : \bar{{\mathcal A}}^{\circ} \to \Delta^{2}\setminus 0$. 
 Since ${\mathcal A}^{\circ} \to \bar{{\mathcal A}}^{\circ}$ is
 a finite morphism, we extend this morhpsim to a finite morphism 
 $\nu : {\mathcal A}' \to \bar{{\mathcal A}}$. 
 The codimension of $\bar{{\mathcal A}}\setminus \bar{{\mathcal A}}^{\circ}$
 is two. By purity of branch loci,
 $\nu$ is \'{e}tale.
 Hence
 $h' : {\mathcal A}' \to \Delta^{2}$ is a smooth abelian fibration.
 By construction, 
 ${\mathcal A}'$ is 
 isomorphic to ${\mathcal A}$ in codimension one.
 Let $A$ be a $h$-ample divisor on ${\mathcal A}$ and $A'$ a
 proper transform on ${\mathcal A}'$. Since every fibre of
 $h'$ is an abelian variety, $A'$ is $h'$-ample. 
 Thus we obtain that ${\mathcal A}$ is isomorphic to ${\mathcal A}'$
 and $h$ is smooth.
 \item Let $A$ be an $f$-ample divisor. We consider the following
       function:
$$
 \lambda (s) := \int_{X_s}\omega \wedge \bar{\omega} A^{\dim S - 2},
$$  
       where $X_s := f^{-1}(s)$ $(s \in S)$.
       Since $f$ is flat, $\lambda (s)$ is a continuous function
       on $S$ by \cite[Corollary 3.2]{fujiki}. 
       Thus $\lambda (s) \equiv 0$ on $S$ and
$$
  \int_{F} \omega \wedge \bar{\omega} A^{\dim S - 2} = 0.
$$
       Since $F$ and $\tilde{F}$ is birational, 
       $j^{*}\omega = 0$ on $\tilde{F}$.
\end{enumerate}
\qqed

\noindent
(2.3) We review basic properties of the mixed hodge structure
      on a simple normal crossing variety.
\begin{lem}\label{restriction}
 Let $X := \sum X_i$ be a simple normal crossing variety.
 Then
$$
 F^{1}H^{1}(X,{\mathbb C}) = 
    \{
     (\alpha_i ) \in \oplus H^{0}(X_i ,\Omega^{1}_{X_i}) |
     \alpha_i |_{X_i \cap X_j } = \alpha_j |_{X_i \cap X_j}
    \}.
$$
\end{lem}
\noindent
{\sc Proof. \quad} 
 Let 
 $$X^{[k]} := 
 \cup_{i_0 < \cdots < i_k} X_{i_0} \cap \cdots \cap X_{i_k} \quad
 \mbox{(disjoint union)}.
 $$
 For an index set $I = \{i_0 , \cdots i_k  \}$, we define 
 an inclusion $\delta^{I}_j$
 $$\delta^{I}_j : X_{i_0} \cap \cdots \cap X_{i_k}
      \to X_{i_1} \cap \cdots \cap X_{i_{j-1}} \cap X_{i_{j+1}} \cap
      \cdots X_{i_k}.$$
 We consider the following spectral sequence \cite[Chapter 4]{GS}:
$$
 E^{p,q}_{1} = H^{q}(X^{[p]}, {\mathbb C}) \Longrightarrow
 E^{p+q} = H^{p+q}(X,{\mathbb C}),
$$
 where $D_2 : E^{p,q}_{1} \to E^{p+1 ,q}_{1}$ is defined by the
 $$
  \bigoplus_{|I| = p} \sum_{j=0}^{p} (-1)^{j}\delta^{I}_{j} .
 $$
 Since
 this spectral sequence degenerates at $E_2$ level (\cite[4.8]{GS}),
 we deduce
$$
{\rm Gr}^{W}_{1}(H^{1}(X,{\mathbb C}))
 = {\rm Ker}(\oplus_{i} H^{1}(X_i ,{\mathbb C}) \stackrel{D_2}{\to}
             \oplus_{i<j} H^{1} (X_i \cap X_j ,{\mathbb C})).
$$
 Moreover $F^{1} \cap W_0 = 0$, $F^{1}H^{1}(X,{\mathbb C}) = 
 F^{1}{\rm Gr}^{W}_{1}(H^{1}(X,{\mathbb C}))$. Thus we obtain
 the assertion of Lemma \ref{restriction} from the definition of $D_2$.
\qqed
\begin{lem}\label{birationa_invariance}
 Let $f: X' \to X$ be a birational morphism 
 between smooth algebraic varieties. Assume 
 the following two conditions:
\begin{enumerate}
 \item There exists a simple normal crossing divisor $Y$ on $X$
       such that $f$ is isomorphic on $X \setminus Y$.
 \item $Y' := (f^{*}Y)_{{\rm red}}$ is a simple normal crossing divisor. 
\end{enumerate}
 Then $F^{1}H^{1}(Y' , {\mathbb C}) \cong F^{1}H^{1}(Y ,{\mathbb C})$.
\end{lem}
\noindent
{\sc Proof. \quad} We consider the following exact sequence of 
 morphisms of Mixed Hodge structures.
$$
     H^{0}(Y' , {\mathbb C}) \stackrel{\alpha}{\to}
     H^{1}(X,{\mathbb C}) \to 
     H^{1}(X' , {\mathbb C})\oplus H^{1}(Y , {\mathbb C}) \to
     H^{1}(Y' , {\mathbb C}) \stackrel{\beta}{\to}
     H^{2}(X,{\mathbb C}) \to
$$
 Note that each morphism has weight $(0,0)$.
 Since $H^{i}(X,{\mathbb C})$ carries pure Hodge structure of weight
 $i$, $\alpha$ and $\beta$ are $0$-map. Moreover 
 $F^{1}H^{1}(X,{\mathbb C}) \cong F^{1}H^{1}(X' ,{\mathbb C})$.
 Thus we deduce that  
 $F^{1}H^{1}(Y' , {\mathbb C}) \cong F^{1}H^{1}(Y ,{\mathbb C})$.
 \qqed
\section{Kulikov model}
\noindent
(3.1)  In this section, we prove the key proposition of the
       proof of Theorem \ref{main}. First we refer the
       the folloing theorem due to Kulikov, 
       Morrison \cite[Classification Theorem I]{morrison}
       and Persson \cite[Proposition 3.3.1]{persson}.
\begin{thm}\label{kulikov} 
 Let $g' : T' \to \Delta$ be a semistable degeneration whose general
 fibre is an abelian surface. Then there exists a semistable
 degeneration $k : {\mathcal K} \to \Delta$ 
 such that $k$ and $g'$ is birational and
 $K_{{\mathcal K}} \sim_{k} 0$. 
 Moreover, exactly one of the following cases occurs:
\begin{enumerate}
 \item ${\mathcal K}_{0}$ is an abelian surface.
 \item ${\mathcal K}_{0}$ consists of a cycle of minimal elliptic ruled
       surfaces, meeting along disjoint sections. 
       Every selfintersection number of double curve is $0$.
 \item ${\mathcal K}_{0}$ consists of a collection of rational surfraces,
       such that the double curves on each component form a 
       cycle of rationall curves; the dual graph $\Gamma$ of 
       $Y^{'}_{0}$ is a triangulation of $S^{1}\times S^{1}$.
\end{enumerate} 
\end{thm}

\noindent
 We call ${\mathcal K}$ a Kulikov model of type I, II or III
 according the case occurs (1), (2) or (3).

\noindent
(3.2) We state the key propositon.
\begin{prop}\label{key}
 Let $f: (X,\omega ) \to S$ be a projective
 complex lagrangian fibration on 4-dimensional symplectic
 manifold
 and $D$ the discriminant locus of $f$.
 If we take a general point $x$  of $D$ and
 an unit disk $\Delta^{1}$ on $S$ such that $\Delta^{1}$ and $D$ intersects
 transversally at $x$
 and  $T := X\times_{S}\Delta^{1}$ is smooth,
 then
\begin{enumerate}
 \item  $t: T \to \Delta^{1}$ 
        is birational to the
        quotient of Kulikov model ${\mathcal K}$ of
        Type I or Type II
        by a cyclic group $G$.
 \item  There exists a nonzero $G$-equivariant element
        of $F^{1}H^{1}({\mathcal K},{\mathbb C})$.
\end{enumerate}
\end{prop}

\noindent
(3.3)
 For the proof of Proposition \ref{key}, we need the following Lemmas.
\begin{lem}\label{nonzero_lemma}
  Let $\nu : Y \to X$ be a birational morphism
  such that $(f\circ \nu)^{*}D$ is a simple normal crossing divisor.
  Then for
 a general point $x$ of 
 the discriminant locus
 $D$ of $f$,  
$$
 F^{1}H^{1}(Y_x , {\mathbb C}) \ne 0
$$
 where $Y_{x} := f^{-1}(x)$.
\end{lem}

\noindent
{\sc Proof. \quad}
 Let $E := ((f\circ \nu)^{*}D)_{{\rm red}}$ and
 $E = \sum E_i$.
 We take an open set $U$ of $S$ which satisfies 
 the following two conditions:
\begin{enumerate}
 \item $D |_{U}$ is a smooth curve.
 \item $f\circ \nu |_{U} : (E |_{f^{-1}(U)})^{[k]} \to D |_{U}$
       is a smooth morphism for every $k$.
\end{enumerate}
 We consider the following exact sequences:
$$
  0 \to {\mathcal F} \to \Omega^{2}_{E_i} \to 
    \Omega^{2}_{E_i /D} \to 0 
$$
$$
 0 \to (f\circ \nu)^{*}\Omega^{2}_{D} \to
    {\mathcal F} 
    \stackrel{\alpha}{\to}
    (f\circ \nu)^{*}\Omega^{1}_{D}\otimes
    \Omega^{1}_{E_i /D} \to 0
$$
 Since $\omega$ is nondegenerate,
 $\nu^{*} \omega \ne 0$ on $E_i$ 
 on a non $\nu$-exceptional divisor $E_i$ .
 By  Lemma \ref{basic_properties} (2),
 the restriction of $\omega$
 on every irreducible component of a fibre of $f$
 is zero. Thus
 $\nu^{*}\omega = 0$ in $ \Omega^{2}_{E_i /D}$.
 On the contrary, $(f\circ \nu )^{*}\Omega^{2}_{D} = 0$,
 we deduce 
$$
 \alpha(\nu^{*}\omega ) \ne 0
$$
 for non $\nu$-exceptional divisor $E_i$.
 Therefore, for a general point $x$ of $D$, 
 $H^{0}(E_{i ,x}, \Omega^{1}_{E_{i , x}}) \ne 0$ where
 $E_{i,x}$ is the fibre of $E_i \to D$ over $x$.
 We denote by $\alpha_i$ the restriction $\alpha(\nu^{*}\omega)$
 to $E_{i , x}$. 
 If $E_{i , x} \cap E_{j , x} \ne \emptyset$, 
 $\alpha_i =\alpha_j$ on $E_{i, x} \cap E_{j , x} $
 by the construction.
 By Lemma \ref{restriction}, we deduce that
 $F^{1}H^{1}(E_{x} , {\mathbb C}) \ne 0$.
\qqed 
\begin{lem}\label{uniqueness_of_Kulikov_model}
  Let $k : {\mathcal K} \to \Delta^{1}$ be a Kulikov model of 
  type I or type II.
  Assume that $k$ is birational to
  a projective abelian fibration $t' :T' \to \Delta^{1}$. Then
\begin{enumerate}
 \item $k$ is a projective morphism.
 \item Every birational map ${\mathcal K} \dasharrow {\mathcal K}$
       is birational morphism.
\end{enumerate}
\end{lem}

\noindent
{\sc Proof. \quad}  
 We may assume that $T'$ is a relatively minimal model over $\Delta^{1}$.
 Then $T'$ and ${\mathcal K}$ is isomorphic in codimension one,
 since $T'$ and ${\mathcal K}$ have only terminal singularities
 and $K_{T'}$ is $t'$-nef.
 Let $A'$ be a $t'$-ample divisor on $T'$ and $A$ a
 proper transform on ${\mathcal K}$. 
 Since $T'$ and ${\mathcal K}$ is isomorhpic in codimension one,
 $A$ is $k$-big.
 If $A$ is $k$-nef, we conclude that $T'$ is isomorhpic to
 ${\mathcal K}$ by relative base point 
 free theorem \cite[Theorem 3-1-2]{KMM}
 Then $k$ is projective,
 ${\mathcal K}$ is the unique relative minimal model 
 and every birational map ${\mathcal K} \dasharrow {\mathcal K}$
 is birational morphism.
 Thus we will prove that $A$ is $k$-nef. Since every big divisor
 on abelian surface is ample, $A$ is $k$-nef if
 ${\mathcal K}$ is
 of type I. In the case that ${\mathcal K}$ is of type II,  
 we investigate the nef cone of each component of the central 
 fibre of ${\mathcal K}$. Let $V$ be a component of the
 central fibre.
 Then $K_V \sim -2e$, where $e$ is
 a double curve. Since $e$ is a section and $e^2 = 0$,
 the nef cone of $V$ is spaned by $e$ and a fibre $l$ of 
 the ruling of $V$.
 We deduce that every effective divisor on $V$ is nef. 
 Therefore $A$ is $k$-nef in the case that
 ${\mathcal K}$ is of type II.
\qqed

\noindent
(3.4) 
{\sc Proof of Proposition \ref{key}. \quad}
 Let $\nu : Y \to X$ be a birational morphism such that
 $(f\circ \nu )^{*}D$ is a simple normal crossing divisor.
 We define $T'' := \nu^{*}T$. If we choose $x$ generally,
 $F^{1}H^{1}(T^{''}_{0}, {\mathbb C}) \ne 0$
 by Lemma \ref{nonzero_lemma}. 
 By Semistable reduction theorem
 \cite[Theorem $11^{*}$]{KKMS}, 
 there exists a generically finite surjective
 morphism $\eta :  T' \to T''$ such that 
 $t' : T' \to \Delta^{1}$ is a semistable degeneration.
 By Theorem \ref{kulikov}, 
 there exists the Kulikov model $k : {\mathcal K} \to \Delta^{1}$ 
 which is birational to $t'$.
 We denote by ${\mathcal K}_{0}$ the central fibre of
 ${\mathcal K}$.
 Then $F^{1}H^{1}({\mathcal K}_0 ,{\mathbb C}) \ne 0$ since
 $F^{1}H^{1}({\mathcal K}_0 ,{\mathbb C})
 \cong F^{1}H^{1}(T^{'}_{0},{\mathbb C})$
 due to Lemma \ref{birationa_invariance} and 
 $F^{1}H^{1}(T^{''}_{0} ,{\mathbb C}) \ne 0$.
 Thus
 ${\mathcal K}$ is of type I or type II. 
 Let 
 $G$ be the galois group of a cyclic extension $K(T')/K(T)$ and
 $g$ a generator of $G$.
 Since $k$ is birational to $t'$, there exists
 a biratinal map $\Phi_{g} : {\mathcal K} \dasharrow {\mathcal K}$ 
 correponding to $g$.
 By Lemma \ref{uniqueness_of_Kulikov_model}, $\Phi_{g}$ is
 a birational morphism and
 $G$ acts on ${\mathcal K}$
 holomorphically. 
 Therefore
 $T$ is birational to the quotient
 ${\mathcal K}/G$.
 We claim that $F^{1}H^{1}({\mathcal K}_{0}, {\mathbb C})$
 carries a nonzero $G$-equivariant element.
 Let $Z$ be a $G$-equivariant resolution 
 of indeterminancy of $T' \dasharrow {\mathcal K}$.
 Then $F^{1}H^{1}(T^{'}_{0}, {\mathbb C}) \cong
       F^{1}H^{1}(Z_0  ,{\mathbb C}) \cong
       F^{1}H^{1}({\mathcal K}_{0}, {\mathbb C})$
 by Lemma \ref{birationa_invariance}.
 Let $\alpha$ be a nonzero element of 
 $F^{1}H^{1}(T^{''}_{0}, {\mathbb C})$. 
 The pull back of $\alpha$ in
 $F^{1}H^{1}(Z_{0},{\mathbb C})$ is
 non zero $G$-equivariant element. Thus
 there exists a nonzero $G$-equivariant element in 
 $F^{1}H^{1}({\mathcal K}_{0}, {\mathbb C})$. 
\qqed 
\section{Classification of Type I degeneration}
\noindent
 (4.1) In this section, we 
 prove the following proposition.
\begin{prop}\label{classification_I}
 Let $t: T \to \Delta^{1}$ be an abelian fibration which is
 birational to the quotient of a Kulikov model ${\mathcal K}$
 of type I by a cyclic group $G$.
 Assume that 
\begin{enumerate}
 \item $T$ is smooth.
 \item $K_T \sim_{t} 0$.
 \item There exists a nonzero $G$-equivariant element of
       $F^{1}H^{1}({\mathcal K}_{0}, {\mathbb C})$.
\end{enumerate}
 Then the representation 
 $\rho : G \to {\rm Aut} H^{1}({\mathcal K}_{0},{\mathbb C})$
 is faithful and
 the central fibre $T_0$ of $T$ satisfies the properties
 of Theorem \ref{main} {\rm (1)}
\end{prop}
\noindent
(4.2) We need the following lemma 
      to prove Proposition \ref{classification_I}. 
\begin{lem}\label{equivariant_fibration}
 The central fibre ${\mathcal K}_{0}$
 admits an $G$-equivariant elliptic fibration
 over an elliptic curve.
\end{lem}
\noindent
{\sc Proof. \quad}
  Since ${\mathcal K}_{0}$ is an abelian surface, it is
  enough to prove that ${\mathcal K}_{0}$ admits
  a $G$-equivariant fibration.
  Let $g$ be a
  generator of $G$. We consider the following morphism:
$$
 H^{1}({\mathcal K}_{0} , {\mathbb C}) \oplus
 H^{1}({\mathcal K}_{0} , {\mathbb C})
 \stackrel{({\rm id}-g^{*})}{\to}
 H^{1}({\mathcal K}_{0} , {\mathbb C}).
$$ 
 Since $G$ is a cyclic group, the kernel of ${\rm id} - g^{*}$
 is $G$-invariant. Moreover this kernel is nonzero by Proposition \ref{key}.
 Therefore $H^{1}({\mathcal K}_{0} ,{\mathbb C})$
 has a $G$-equivariant sub Hodge structure and we 
 conclude that ${\mathcal K}_{0}$ admits a
 $G$-equivariant fibration. 
\qqed

\noindent
(4.3) {\sc Proof of Proposition \ref{equivariant_fibration}. \quad}
 We will construct
 a suitable resolution $Z$ of ${\mathcal K}/G$ and the unique relative
 minimal model $W$ of $Z$ over $\Delta^{1}$. 
 By Lemma \ref{equivariant_fibration}, 
 there exists a $G$-equivariant elliptic fibration
 on the central fibre ${\mathcal K}_{0}$ of ${\mathcal K}$.
 We denote this fibration by  $\pi : {\mathcal K}_{0} \to C$.
 By construction, the  action of $G$ on $C$ is translation.
 Let $H$ be the kernel of
 the representation 
 $\rho : G \to {\rm Aut}H^{1}({\mathcal K}_0 , {\mathbb C})$.
 Since the action of $H$ on ${\mathcal K}_0$ is 
 a translation, 
 ${\mathcal K}/H$ is smooth. 
 It is enough to consider the action of $G/H$ on
 ${\mathcal K}/H$ for the investigation of 
 the singularities of ${\mathcal K}/G$. 
 If the action of $G/H$ on $C/H$ is translation,
 then ${\mathcal K}/G$ is smooth. Moreover ${\mathcal K}/G$
 is the unique relative minimal model over $\Delta^{1}$ since
 it has no rational curves.
 On the contrary, $T$ is a relative minimal model over $\Delta^{1}$,
 $T \cong {\mathcal K}/G$. 
 By construction,
 the central fibre of the quotient ${\mathcal K}/G$
 is an hyperelliptic surface. 
 Since
 every hyperelliptic surface is the \'{e}tale quotient
 of the product of elliptic curves, 
 $T_0 \cong ({\mathcal K}/G)_{0}$ satisfies the property 
 of Theorem \ref{main}.
 We claim that the representation
 $\rho : G \to {\rm Aut}H^{1}({\mathcal K}_0 , {\mathbb C})$
 is faithful.
 If $H$ is not trivial,
 then ${\mathcal K}$ is not $G$-invariant and 
 $K_{{\mathcal K}/G} \not\sim 0$. Since 
 ${\mathcal K}/G \cong T$ and $K_T \sim 0$, $H$ is trivial.
 In the following, we assume that 
 the action of $G/H$ on $C/H$ is trivial. 
 Since $\pi : {\mathcal K}_{0} \to C/H$ is $G/H$-equivariant,
 the singularites of ${\mathcal K}/G$ consists of
 several copies the product of a surface quotient singularity 
 and an elliptic curve. 
 The list of surface quotient
 singularities which occur above
 is found in \cite[Table 5 (p157)]{BPV}.
 We construct the relative minimal
 resolution $Z$ of ${\mathcal K}/G$ by
 the minimal resolution of surface quotient singularities.
 If the singularities of ${\mathcal K}/G$
 consists of the product of Du Val singularities and
 an elliptic curve only, then
 $Z$ is a relative minimal model over $\Delta^{1}$ and
 we put $W=Z$.
 In other cases, we obtain a relative minimal model $W$
 after birational contractions of $Z$ (cf. \cite[pp 156--158]{BPV}).
 In both cases,
 $W$ is the unique minimal model
 by the similar argument in Lemma \ref{uniqueness_of_Kulikov_model}.
 Since $W$ is birational to $T$ and
 $T$ is a relative minimal model over $\Delta^{1}$,
 $T \cong W$.
 By construction, the central fibres  $W_{0}$ admit a
 fibration over $C/H$. Note that
 the fibre of $W_0 \to C/H$ is 
 a Kodaira singular fiber of type  $I^{*}_{0}$, $II$,
 $II^{*}$, $III$, $III^{*}$, $IV$ or $IV^{*}$.
 Since ${\rm Sing}W_{0}$ forms multi sections of
 $W_0 \to C/H$, there exists an \'{e}tale
 finite cover $\tilde{C} \to C/H$ and the base change
 $W_{0}\times_{C/H}\tilde{C}$ is isomorphic to
 the product of a Kodaira singular fibre and an elliptic curve.
 Finally, we prove that the representation 
 $\rho : G \to {\rm Aut} H^{1}({\mathcal K}_{0},{\mathbb C})$
 is faithful. We derive a contradiction assuming
 that $H$ is not trivial.
 If $W=Z$, then there exists
 a morphism $\eta : W \to {\mathcal K}/G$
 such that $\eta^{*}K_{{\mathcal K}/G} \sim K_W$. Since the 
 action of $H$ on ${\mathcal K}_{0}$ is translation, $K_{{\mathcal K}}$
 is not $G$-equivariant and $K_{{\mathcal K}/G} \not\sim 0$.
 However, $K_W \sim \eta^{*}K_{{\mathcal K}/G}$ and $K_W \sim 0$,
 that is a contradiction. If $W \not\cong Z$
 we consider the base change 
 $Z_1 := Z\times_{{\mathcal K}/G}{\mathcal K}/H$.
 Since $Z$ is obtained by 
 blowing up along singular locus
 of ${\mathcal K}/G$ (cf. \cite[pp 158]{BPV})
 and the singular locus of ${\mathcal K}/G$ consists of
 elliptic curves,
 $Z_1$ is smooth. Let $\eta_1 : Z_1 \to W$, 
 $\eta_2 : Z_1 \to {\mathcal K}/H$ and 
 $k/H : {\mathcal K}/H \to \Delta^{1}$. 
 We denote by $F$ the central fibre of $k/H$ 
 with reduced structure.
 Note that $(k/H)^{*}(0) = mF$, where $m$
 is the order of $H$.
 Then 
$$
  \eta^{*}_{1}K_W  \sim \eta_{2}^{*} K_{{\mathcal K}/H} - \eta_{2}^{*}F.
$$
 By adjunction formula,
 $K_{{\mathcal K}/H} \sim (m-1)F$.
 Thus $K_W \not\sim 0$.
 However, this is a contradiction because $K_W \sim 0$.
\qqed

\section{Classification of Type II degeneration}
\noindent
(5.1) In this section, we prove the following proposition and
      Theorem \ref{main}.
\begin{prop}\label{classification_II}
 Let $t: T \to \Delta^{1}$ be an abelian fibration which is
 birational to the quotient of a Kulikov model ${\mathcal K}$
 of type II by a cyclic group $G$.
 Assume that 
\begin{enumerate}
 \item $T$ is smooth.
 \item $K_T \sim_{t} 0$.
 \item There exists a nonzero $G$-equivariant element of
       $F^{1}H^{1}({\mathcal K}_{0}, {\mathbb C})$.
\end{enumerate}
 Then the representation $\rho : G \to 
 {\rm Aut} H^{1}({\mathcal K}_{0},{\mathbb C})$ is faithful and
 the central fibre $T_0$ of $T$ satisfies the
 properties of Theorem \ref{main} {\rm (2)}.
\end{prop}
\noindent
(5.2) For the proof of Proposition \ref{classification_II},
      we investigate the action of $G$ on 
      the central fibre of ${\mathcal K}$.
\begin{lem}\label{cyclic_case}
 Let $g$ be a generator of $G$ and $m$ the smallest positive
 interger such that
 every component is stable under the action of $g^m$.
 We denote by $H$ the subgroup of $G$ 
 generated by $g^m$. 
 Then
\begin{enumerate}
 \item The representation 
       $\sigma : G \to {\rm Aut}F^{1}H^{1}({\mathcal K}_{0}, {\mathbb C})$
       is trivial. 
 \item The action of $H$ is free and the central
       fibre of the quotient ${\mathcal K}/H$
       is a cycle
       of mininal elliptic ruled surfaces.
\end{enumerate}
\end{lem}
\noindent
{\sc Proof. \quad} 
\begin{enumerate}
 \item By Proposition \ref{birationa_invariance},
       there exists a $G$-equivariant element in 
       $F^{1}H^{1}({\mathcal K}_{0}, {\mathbb C})$. 
       Since $\dim F^{1}H^{1}({\mathcal K}_{0}, {\mathbb C}) = 1$,
       every element of $F^{1}H^{1}({\mathcal K}_{0}, {\mathbb C})$
       is $G$-invariant.
 \item From the assumption 
       there exists an action of $g^m$ on each component
       of the central fibre of ${\mathcal K}$. 
       Let $V$ be a component of the central fibre
       and $\pi : V \to C$ a ruling.
       Since every fibre of $\pi$ is ${\mathbb P}^{1}$ 
       and $C$ is an elliptic curve,
       $g^m$ maps a fibre of $\pi$ to a fibre of $\pi$, that
       is, $\pi$ is $g^m$-equivariant.
 From Lemma \ref{cyclic_case} (1) and
 Lemma \ref{restriction}, 
 holomorphic one forms on $V$
 is invariant under the action of $g^m$.
 Thus, the action of $g^m$
 on $C$ is translation.
 Therefore $V/H$ is  
 a minimal elliptic ruled surface. From the assumption
 that each component is stable under the action of $g^m$,
 the central fibre of the quotient
 ${\mathcal K}/H$ 
 is a cycle of 
 minimal elliptic ruled surfaces. \qqed
\end{enumerate}

\noindent
(5.3)  {\sc Proof of Proposition \ref{classification_II}. \quad}
 From Lemma \ref{cyclic_case}, ${\mathcal K}/H $ is smooth and
 the central fibre of ${\mathcal K}/H$ is a cycle of minimal elliptic ruled 
 surfaces. 
 Let $\Gamma$ be the dual graph of the central fibre of ${\mathcal K}$
 and $g$ a generator of $G$.
 Considering ${\mathcal K}/H$ in stead
 of ${\mathcal K}$, we may assume that the action of $g^{m}$ 
 is trivial if the action of $g^{m}$ on $\Gamma$ is trivial.

\noindent
(5.3.1) 
 If the action of $G$ is free,
 ${\mathcal K}/G$ is smooth and
 this is a relative minimal model over $\Delta^{1}$.
 Every component of the central fibre of ${\mathcal K}/G$
 is the minimal elliptic ruled surface $V$ which has
 a section $e$ such that $K_V \sim -2e$ and $e^2 = 0$.
 Thus we show that
 $T \cong {\mathcal K}/G$ by
 similar argument as in the proof of
 Lemma \ref{uniqueness_of_Kulikov_model}.
 We claim that the representation
 $\rho : G \to 
       {\rm Aut} H^{1}({\mathcal K}_{0},{\mathbb C})$
 is faithful. If ${\rm Ker}\rho$ is not trivial,
 then $K_{{\mathcal K}}$ is not $G$-equivariant and
 $K_{{\mathcal K}/G} \not\sim 0$. However $K_T \sim 0$,
 we obtain ${\rm Ker}\rho$ is trivial.
 Since $\Gamma$ is a Dynkin diagram of type $\tilde{A_n}$
 and $G$ is a cyclic group,
 the action of $G$ on $\Gamma$ is either rotation or reflection.
\begin{enumerate}
 \item If the action of $G$ on $\Gamma$ is rotation,
       the central fibre of
       ${\mathcal K}/G$ is a cycle of minimal elliptic ruled surfaces.
       Each double curve is a section of a minimal elliptic ruled
       surface.
 \item If the action of $G$ on $\Gamma$ is reflection,
       the central fibre $({\mathcal K}/G)_{0}$ of
       ${\mathcal K}/G$ is a chain of minimal elliptic ruled surfaces.
       Let $V$ be the edge component of the central fibre
       and $V'$ the component such that $V \cap V' \ne \emptyset$.
       Note that 
$({\mathcal K}/G)_{0} = 2V + 2V' + \mbox{(other components)}$.
       By adjunction formula,
$$
       K_V \equiv -V'|_{V}.
$$
       Therefore, the double curve $V\cap V'$ is a
       bisection of the ruling of $V$. Every other double curve
       is a section.
\end{enumerate}

\noindent
(5.3.2) If the action of $g$ is not free,
 we need the following lemma.
\begin{lem}
 If the action of $G$ has fixed points,
 then the
 action of $G$ on $\Gamma$ is reflection and it
 preserves
 two vertices. 
 The fixed locus of the
 central fibre consists of sections or bisections of the
 ruling of components corresponding to
 the fixed vertices.
%
%
\end{lem}
\noindent
 Assuming this Lemma, the central fibre of the quotient ${\mathcal K}/G$
 is a chain of minimal elliptic ruled surfaces. The singularities
 of ${\mathcal K}/G$ consists of several
 copies of the product of $A_1$ singlarity and
 an elliptic curve. Thus the unique relative minimal
 model $W$ over $\Delta^{1}$ is obtained by blowing up along singlar locus.
 Since $T$ is a relative minimal model over $\Delta^{1}$, 
 $W \cong T$. The dual graph of the central fibre of $W$
 is $A_n$, $D_n$ or $\tilde{D_n}$. The double curve 
 on the edge component is a bisection or a section.
 Every other double curve is section.
 We claim that the representation
 $\rho : G \to {\rm Aut} H^{1}({\mathcal K}_{0},{\mathbb C})$
 is faithful. 
 If $H$ is not trivial, $K_{{\mathcal K}}$ is not $G$-equivariant.
 Since $K_W$ is the pull back of
 $K_{{\mathcal K}/G}$ is crepant, $K_W \not\sim 0$ if
 $H$ is not trivial. However, $K_T \sim 0$, that is a contradiction.

\noindent
(5.4) {\sc Proof of Lemma \ref{cyclic_case}. \quad}
 If the action of $G$ on $\Gamma$ is rotation, there exists no
 fixed points. Thus the action of $G$ on $\Gamma$ is reflection. 
 We derive
 the contradiction assuming that $G$ fixes one of edges of $\Gamma$.
 Let $C$ be the elliptic curve corresponding to the edge which
 is fixed by $G$. From Lemma \ref{restriction} and
 Lemma \ref{cyclic_case},
 the action of $G$ on $C$ 
 preserves holomorphic one form on $C$. Therefore
 $C$ is fixed locus of the action of $G$. 
 The singularities of the quotient ${\mathcal K}/G$ consist of
 several copies of the product of $A_1$ singularity and
 an elliptic curve. Let
 $w: W \to {\mathcal K}/G$ be the blowing up 
 along $C$. 
 We denote by $V_i$
 each components of $W_0$. Let $V_0$ be
 the exceptional divisor coming from the blowing up 
 along $C$. Since 
 the central fibre $W_0$ of $W$ is a chain of minimal
 ellitic ruled surfaces,
 there exists componets 
 $V_1$ and $V_2$
 such that $V_0 \cap V_1 \ne \emptyset$,
 $V_2 \cap V_1 \ne \emptyset$ and
 $V_i \cap V_j = \emptyset$ $(i = 0,1,2, j\ne 0,1,2)$.
 Then 
 $W_0 = V_0 + 2V_1 + 2V_2 + \mbox{(Other component)}$.
 Since $W$ is smooth along $V_1$,
$$
 K_{V_1} \equiv K_W + V_1 |_{V_1} \equiv 
 (-\frac{1}{2}V_0 -V_2 )|_{V_1}.
$$
 by adjunction formula.
 Let $l$ be a fibre of ruling of $V_1$. Then
$$
 K_l \equiv K_{V_1} + l |_{l} \equiv (-\frac{1}{2}V_0 - V_2).l .
$$
 Since every double curve of $W_0$ is a section, 
 $\deg K_l = -3/2$. However this is a contradiction
 because $l \cong{\mathbb P}^{1}$. Therefore $G$ fixes
 two vertices. Let $V$ be one of the component correponding to
 the fixed vertices and $\pi : V \to C$ the 
 ruling of $V$. Since $C$ is an elliptic curve and every
 fibre $\pi$ is ${\mathbb P}^{1}$, 
 $\pi$ is $G$-equivariant.
 By Lemma \ref{cyclic_case} (1), the action of $G$ on $V$ preserves
 one form on $V$.
 Since the action of $G$ is not free,
 $G$ acts on $C$ tirivially. 
 There exist two fixed points
 on each fibre of the ruling of $V$. 
 Thus we obtain
 the rest of assertion of Lemma \ref{cyclic_case}.
\qqed

\noindent
 The proof of Proposition \ref{classification_II} is completed.
\qqed

\noindent
(5.5) {\sc Proof of Theorem \ref{main}. \quad}
 We take a general point $x$ of the
 discriminant locus of $f$ and an unit disk $x \in \Delta^{1}$ such 
 that $T := X\times_{S} \Delta^{1}$ is smooth.
 By Proposition \ref{key}, the abelian fibraton $T \to \Delta^{1}$
 satisfies assumptions of Proposition 
 \ref{classification_I} or \ref{classification_II}.
 Then $T_0$ satisfies the assertions of Theorem \ref{main}
 by Proposition \ref{classification_I} and \ref{classification_II}.
\qed 
 
\newpage

\end{document}